\documentclass{esub2acm}
 \usepackage{amssymb}

 \newcommand{\be}{\begin{equation}}
 \newcommand{\ee}{\end{equation}}
 \newcommand{\ea}{\end{array}}
 
 \newcommand{\Ai}{\mbox{Ai}}
 \newcommand{\Bi}{\mbox{Bi}}

 \def\Frac#1#2{\frac{\displaystyle{#1}}{\displaystyle{#2}}}
 
\begin{document}

 \title{Computing solutions of the modified 
Bessel differential equation for imaginary orders
and positive arguments}
 \author{Amparo Gil \footnote{Present address: Departamento de Matem\'aticas, 
 Estad\'{\i}stica y Computaci\'on. 
 U. de Cantabria, 39005-Santander, Spain}\\
         Departamento de Matem\'aticas, U. Aut\'onoma de Madrid, 
         28049-Madrid, Spain\\
         \\
     Javier Segura\\
         Departamento de Matem\'aticas, Estad\'{\i}stica y 
 Computaci\'on. 
 U. de Cantabria, 39005-Santander, Spain\\   
     \and
     Nico M. Temme\\
     CWI, P.O. Box 94079, 1090 GB Amsterdam, The Netherlands}

 \begin{abstract}
We describe a variety of methods to compute the functions
$K_{ia}(x)$, $L_{ia}(x)$ and their derivatives for real $a$ 
and positive $x$. These functions are numerically 
satisfactory independent solutions
 of the differential equation  $x^2 w'' +x w' +(a^2 -x^2)w=0$.
 In the accompanying paper \cite{Gil0}, we describe the 
implementation of these methods in Fortran 77 codes. 
 \end{abstract}
 \begin{bottomstuff}

 \end{bottomstuff}
 \markboth{Amparo Gil, Javier Segura and Nico M. Temme}
      {Modified Bessel functions of imaginary orders}
 \maketitle
 \section{Introduction}

In previous publications \cite{Gil1, Gil2}, methods to compute the modified
 Bessel function $K_{ia}(x)$ for positive $x$ were described. We complete
 here this analysis by describing analogous methods for the computation of
 the function $L_{ia}(x)$. With this, methods for the reliable computation
 of a pair of linearly independent numerically satisfactory solutions become
 available which find their implementation in the 
 accompanying paper \cite{Gil0}. Methods to compute their derivatives
are also provided. 

 The functions $K_{ia}(x)$ and $L_{ia}(x)$ are solutions of the 
 modified Bessel equation
 for imaginary orders

\begin{equation}
x^2 w'' +x w' +(a^2 -x^2)w=0.
\label{ODE}
\end{equation}
\noindent
The function $K_{ia}(x)$ finds application in a number of problems of
physics and applied mathematics \cite{Gil1}. The function $L_{ia}(x)$ 
is a real valued numerically satisfactory companion to $K_{ia}(x)$ in
the sense described
in \cite{Olv}, pp. 154--155.

In terms of the modified Bessel function of the first kind $I_{\nu}(x)$, 
the solutions are defined as:

\begin{equation}
K_{i a}(x)=\Frac{\pi}{2i\sinh(\pi a)}
\left[I_{-ia}(x)-I_{ia}(x)\right]\,,\,\,
L_{ia}(x)=\Frac{1}{2}\left[I_{-ia}(x)+I_{ia}(x)\right],
\label{definition}
\end{equation}
\noindent with Wronskian $W\left[K_{ia}(x),L_{ia}(x)\right]=1/x$. 

Both $K_{ia}(x)$
and $L_{ia}(x)$ are real solutions for real $x>0$ and $a\in {\mathbb R}$. 
Because they are even functions of $a$, in the sequel we will consider
$a\ge 0$, although this restriction is not present in the code.

In Section 2, we describe the different methods of computation considered, namely:
series expansions, asymptotic expansions for large $x$, Airy-type uniform 
asymptotic expansions, non-oscillating integral representations (including a discussion
on the quadrature rule) and a continued fraction method. We avoid duplicating information
already given in previous papers; in particular, the references \cite{Gil1,Gil2} provide
information required for building the algorithms of the accompanying paper \cite{Gil0}. 
A few misprints in \cite{Gil1} are corrected.

In Section 3, we include a discussion on the dominant asymptotic behaviour of the functions.
These exponential dominant factors can be taken out, leading to scaled functions which
can be computed in a much wider range. The algorithm described in the accompanying paper 
\cite{Gil0} offers the possibility of computing scaled and unscaled functions.

 \section{Methods of computation}
\label{metodos}

In \cite{Gil1,Gil2} methods are described to compute the $K_{ia}(x)$
for different regions in the $(x,a)$ plane. In particular, we considered 
series expansions \cite{Tem1},
asymptotic expansions for large $x$ (\cite{Abr}, Eq. 9.7.2), 
uniform asymptotic expansions 
for $a\simeq x$ (\cite{Bal,Dun} and \cite{Olv}, pg. 425). Also,  
non-oscillating integral representations \cite{Tem2,Gil1} are available.
 Similar techniques are available 
for the computation of $L_{ia}(x)$ and the derivatives 
$K'_{ia}(x)$, $L'_{ia}(x)$. In addition,
a continued fraction method can be applied for the computation of $K_{ia} (x)$ and $K'_{ia}(x)$. 
Those techniques generally give at least two alternative methods for computing the functions
in the $(x,a)$ plane for moderate  values of $x$ and $a$; therefore, we
can always compare different methods of computation for testing 
their accuracy. The selection of one or
another method of computation in a given region will depend on the range of applicability of each method and
its efficiency.

We now describe the different methods of computation which are used
in the programs.

\subsection{Series expansions}

Series expansions can be built which properly describe the solutions near the singular
point ($x=0$) of the defining differential equation (\ref{ODE}). The idea, as described 
in \cite{Tem1,Gil1}, is to substitute the Maclaurin series for $I_{\nu}(x)$ 
(\cite{Abr}, Eq. 9.6.10) in Eqs. (\ref{definition}). The following expansions are obtained

\begin{equation}
\begin{array}{lll}
K_{ia}(x)=\Frac{1}{n(a)}\displaystyle\sum_{k=0}^{\infty}f_k c_k &, & K'_{ia}(x)=\Frac{1}{n(a)}
\Frac{2}{x}\displaystyle\sum_{k=0}^{\infty}
\left[k f_k-\Frac{r_k}{2}\right]c_k\\
L_{ia}(x)=n(a)
\displaystyle\sum_{k=0}^{\infty}r_k c_k &, & L'_{ia}(x)=  n(a)
\Frac{2}{x}\displaystyle\sum_{k=0}^{\infty}
\left[k r_k+a^2 \Frac{f_k}{2}\right]c_k ,
\end{array}
\label{series}
\end{equation}
\noindent
where 

\begin{equation}
n(a)=e^{\pi a/2}\sqrt{\Frac{1-e^{-2\pi a}}{2\pi a}}\,,\quad
c_k=(x/2)^{2k}/k!
\end{equation}
and \cite{Dun}
\begin{equation}
\begin{array}{l}
f_k=\Frac{\sin(\phi_{a,k}-a\ln (x/2))}{(a^2 (1+a^2)...(k^2+a^2))^{1/2}} ,\\
\\
f_k/r_k=\Frac{1}{a}\tan(\phi_{a,k}-a\ln (x/2)), \mbox{ with } 
\phi_{a,k}=\arg\,(\Gamma (1+k+ia)) .
\label{explicit}
\end{array}
\end{equation}
\noindent 
The coefficients $f_k$ and $r_k$ differ from those in \cite{Gil1} by a constant factor (for fixed $a$).
The new normalization shows explicitly (Eqs. (\ref{series})) the dominant exponential  behaviour $n(a)$ and $1/n(a)$ as 
$a\rightarrow\infty$ ($\sim e^{\pm \pi a/2}$). 

An efficient method to compute the coefficients was described in \cite{Gil1,Tem1}; this method is
 based on the fact that
both $f_k$ and $r_k$ satisfy the three-term recurrence relation
\begin{equation}
(k^2+a^2)r_k-(2k-1)r_{k-1}+r_{k-2}=0 .
\end{equation}
Perron's theorem is inconclusive with respect to the existence of minimal solutions for this recurrence relation;
anyhow, 
the second equation in (\ref{explicit}) confirms that neither $f_k$ nor $r_k$ are minimal solutions.
Therefore, forward recursion will be numerically stable. Starting values can be computed taking into account
that $\arg\,\Gamma (1+ia)=\sigma_0 (a)$, where $\sigma_0 (a)$ is the Coulomb phase shift, for which Chebyshev
expansions are available for double precision \cite{Cod}. Namely, we have:

\begin{equation} \begin{array}{l}
r_0=\cos[\sigma_0 (a)-a\ln (x/2) ],\\ \\
r_1=\Frac{1}{1+a^2}\left\{\cos [\sigma_0 (a)-a\ln (x/2)]- 
a \sin[\sigma_0 (a)-a\ln (x/2)
]\right\}, \end{array} \end{equation}

and
\begin{equation} \begin{array}{l}
f_0=\Frac{1}{a}\sin[\sigma_0 (a)-a\ln (x/2) ],\\ \\
f_1=\Frac{1}{a(1+a^2)}\left\{\sin [\sigma_0 (a)-a\ln (x/2) ]
+ a \cos[\sigma_0 (a)-a\ln (x/2)]\right\}\,. 
\end{array} \end{equation}

These formulas correct two misprints in \cite{Gil1} (Eqs. 12 and 13). 

Series can be used for $x/a$ small. See \cite{Gil0}, Section 2 \label{canv}.

\subsection{Asymptotic expansions for large $x$}
\label{asyx}

Asymptotic expansions for large $x$ are available from the known
asymptotic expansion of $I_{\nu}(x)$ (\cite{Abr}, Eq. 9.7.1):
\begin{equation}
\begin{array}{l}
K_{ia}(x)= 
\left(\Frac{\pi}{2x}\right)^{1/2}e^{-x}
\left\{\displaystyle\sum_{k=0}^{n-1}\Frac{(ia,k)}{(2x)^k}+\gamma_n
\right\},\\
\\
L_{ia}(x)= \Frac{1}{\sqrt{2\pi x}}e^{x}
\left\{\displaystyle\sum_{k=0}^{n-1}(-1)^k \Frac{(ia,k)}{(2x)^k}
+\delta_n\right\},
\end{array}
\label{largex}
\end{equation}
\noindent
where $(ia,m)$ is the Hankel symbol, which satisfies 
\begin{equation}
(ia,k+1)=-\Frac{\left(k+\frac12 \right)^2 +a^2}{k+1}(ia,k)\,,\,\, (ia,0)=1 .
\end{equation}

Bounds for the error terms ($\gamma_n$, $\delta_n$) 
can be found in \cite{Olv}, Pg. 269, Ex. 13.2.

As discussed in \cite{Gil1} the numerical performance of the 
asymptotic expansion for $K_{ia}(x)$ is of more
restricted applicability than for the case of the evaluation of $K_{\nu}(x)$ for real
$\nu$. Furthermore, the
continued fraction method described 
in \cite{Gil1} covers the region where this 
expansion is of
numerical interest. For this reason, the continued fraction method is the preferred algorithm for the
evaluation of $K_{ia}(x)$ and $K^{\prime}_{ia}(x)$ for moderate
values of $a$. On the other hand the asymptotic expansion 
for $L_{ia}(x)$ turns out to be
 accurate in a wider region,
which is a fortunate situation given that the continued fraction method 
is not available in this case. See \cite{Gil0}, Section 2, for further details.

Asymptotic expansions for the derivatives are also available by differentiating Eqs. (\ref{largex}).

\subsection{Airy-type uniform asymptotic expansions}
\label{atuae}

The Airy-type asymptotic expansions for $K_{ia}(x)$ can be found 
in \cite{Bal,Dun} and \cite{Olv} pg. 425; the analogous expansions 
for $L_{ia}(x)$ \cite{Dun} 
and $K^{\prime}_{ia}(x)$ are also available \cite{Bal}, while the expansion
for $L^{\prime}_{ia}(x)$ can be derived in the same way. We summarize here
the main features needed for the computation through these expansions, 
neglecting the error terms.
Further details can be found in \cite{Bal,Dun,Olv} and \cite{Gil2,Tem3}.

  The expansion for $K_{ia}(x)$ and $L_{ia}(x)$ in terms of Airy functions ($\Ai(z)$, $\Bi (z)$ and their derivatives)
reads

\begin{equation}
\begin{array}{l}
K_{ia}(a z)=\Frac{\pi e^{-a\pi/2} \phi(\zeta)}{a^{1/3}}
\left[\Ai(-a^{2/3}\zeta) F_a (\zeta) + \Frac{1}{a^{4/3}}\Ai^{\prime}(-a^{2/3}\zeta) G_a (\zeta)\right] ,\\
\\
L_{ia}(a z)=\Frac{e^{a\pi/2} \phi(\zeta)}{2a^{1/3}}
\left[\Bi(-a^{2/3}\zeta) F_a (\zeta) +  \Frac{1}{a^{4/3}} \Bi^{\prime}(-a^{2/3}\zeta)  G_a (\zeta)\right] ,
\end{array}
\label{unifun}
\end{equation}

\noindent
where 
\begin{equation}
 F_a (\zeta) \sim\displaystyle\sum_{s=0}^{\infty} (-)^s\Frac{a_s(\zeta)}{a^{2s}}\,,\,\,
 G_a (\zeta) \sim\displaystyle\sum_{s=0}^{\infty} (-)^s\Frac{b_s(\zeta)}{a^{2s}} ,
\label{FaGa}
\end{equation}
\noindent as $a\rightarrow \infty$ uniformly with respect to $z\in[0,\infty)$. Error bounds for the asymptotic
expansions of the $K_{ia}(x)$ and $L_{ia}(x)$ are given in \cite{Dun}.

The quantity $\zeta$ is given by

\begin{equation}
\begin{array}{ll}
   \Frac{2}{3}\zeta^{3/2}&=\log\Frac{1+\displaystyle\sqrt{1-z^2}  }{z}-\displaystyle\sqrt{1-z^2},\,\,0\le z\le 1  ,   \\
& \\
 \Frac{2}{3}(-\zeta)^{3/2}&=\displaystyle\sqrt{z^2-1}-\arccos\Frac{1}{z},\,\,z\ge 1 ,\\
\end{array}
\label{zeta}
\end{equation}

and

\begin{equation}
\phi(\zeta)=\left(\Frac{4\zeta}{1-z^2}\right)^{1/4},\,\,\phi(0)=2^{1/3} .
\end{equation}

Of course, it is crucial to compute accurately Eqs. (\ref{zeta}) for small $\zeta$. For this, series expansions around
$z=1$ can be used.

  The evaluation of the coefficients near the turning point $z=1$ (which is our region of interest) can be performed via
 Maclaurin series expansions of the quantities $\phi$, $a_s$ and $b_s$(\cite{Tem3}) in terms of the 
variable $\eta=2^{-1/3}\zeta$ (see \cite{Gil2} and \cite{Tem3} 
for further details).

Asymptotic expansions for the derivatives can be found by differentiating Eqs. (\ref{unifun}). In this way:

\begin{equation}
\begin{array}{l}
K'_{ia}(a z)=2\Frac{\pi e^{-a\pi/2}}{za^{2/3}\phi (\zeta)}
\left[\Ai^{\prime}(-a^{2/3}\zeta) P_a (\zeta)+\Frac{1}{a^{2/3}}\Ai(-a^{2/3}\zeta) Q_a (\zeta)\right] ,\\
\\
L'_{ia}(a z)=\Frac{e^{a\pi/2}}{z a^{2/3}\phi(\zeta)}
\left[\Bi^{\prime}(-a^{2/3}\zeta) P_a (\zeta)+\Frac{1}{a^{2/3}}\Bi(-a^{2/3}\zeta) Q_a (\zeta)\right] ,
\end{array}
\label{unifun2}
\end{equation}

\noindent where $P_a(\zeta)$ and $Q_a (\zeta)$ can be written in terms of $F_a (\zeta)$, $G_a (\zeta)$ and their
derivatives and they have asymptotic expansions

\begin{equation}
 P_a (\zeta) \sim\displaystyle\sum_{s=0}^{\infty} (-)^s\Frac{c_s(\zeta)}{a^{2s}}\,,\,\,
 Q_a (\zeta) \sim\displaystyle\sum_{s=0}^{\infty} (-)^s\Frac{d_s(\zeta)}{a^{2s}} ,
\label{PaQa}
\end{equation}

\noindent
 whose coefficients can be obtained from the computed coefficients $a_s$ and $b_s$ (in Taylor series
around $\zeta =0$) through the relations:

\begin{equation}
\begin{array}{l}
c_s (\zeta )=a_s (\zeta)+\chi (\zeta )b_{s-1}(\zeta )+b'_{s-1} (\zeta) ,\\
d_s (\zeta )=-\chi (\zeta )a_s (\zeta)-a'_s (\zeta)-\zeta b_s (\zeta)  ,
\end{array}
\label{ces}
\end{equation}

\noindent
where
\begin{equation}
\chi (\zeta ) =\phi' (\zeta)/\phi(\zeta),
\label{chi}
\end{equation}
The prime in Eqs. (\ref{ces}) and (\ref{chi}) denotes the derivative with respect to $\zeta$. Using Eqs. (\ref{ces}) the
coefficients $c_s$ and $d_s$ can be computed from the coefficients $a_s$ and $b_s$. Details on the evaluation of 
$a_s$ and $b_s$ are given in \cite{Gil2}, where an explicit Maple algorithm is given for the computation of $a_s$ and $b_s$
for $s=0,1,2,3$.

By computing the Wronskian relation for the modified Bessel functions and using the Wronskian for Bessel functions,
it is easy to derive the relation
\begin{equation}
F_a (\zeta) P_a (\zeta) -\Frac{1}{a^2}G_a (\zeta) Q_a (\zeta)=1
\end{equation}
\noindent
which is a useful relation for checking the correctness of the approximations for the coefficients in the asymptotic
expansions.

An algorithm to compute Airy functions of a real variable is needed for the
computation of these asymptotic expansions. In the routines \cite{Gil0} we use Algorithm 819 \cite{Gil3}.

These Airy-type asymptotic expansion are applied in \cite{Gil0} in a broad region around the turning point line $a=x$.

\subsection{Non-oscillating integral representations}
\label{noir}

Paths of steepest descent for integral representations of the modified Bessel functions of imaginary
orders and their derivatives are given in \cite{Tem2}. Apart from their application in asymptotics \cite{Memo}, 
these integrals are useful for building 
numerically stable (non-oscillating) integral representations 
for $K_{ia}(x)$ and 
$K_{ia}'(x)$, as described in \cite{Gil1}.
We complete here the analysis in \cite{Gil1} by providing analogous expressions for the computation of $L_{ia}(x)$ and
$L'_{ia}(x)$. Additionally, we study further transformations of the integrals which enable us
to obtain integral expressions suitable for computation by means of 
the trapezoidal rule.

\subsubsection{Monotonic case ($x>a$)}

We have the following integral representations
in the monotonic region \cite{Gil1}
\begin{equation}
\begin{array}{l}
K_{ia}(x)=e^{-\lambda}\displaystyle\int_0^{\infty}e^{-x\Phi (\tau)}d\tau\\
\\
K^{\prime}_{ia}(x)=-e^{-\lambda}\,\displaystyle\int_{0}^{\infty}\left[
\cos\theta+\Frac{\cosh\,\tau-1+2\sin^2\frac12(\theta-\sigma)}{\cos
\,\sigma}\, \right]\, e^{-x \Phi(\tau)}\, d \tau,
\end{array}
\label{kmon}
\end{equation}
\noindent 
where 
\begin{equation}
\lambda=x\cos\theta +a\theta\,,\,\,
a=x\sin\theta \,,\,\,\sin
\sigma=\left(\sin\theta \Frac{\tau}{\sinh \tau}\right)
\label{definiciones}
\end{equation}
\noindent and $\theta \in [0,\pi /2)$, $\sigma\in (0,\theta]$. 
The dominant exponential term ($e^{-\lambda}$)
 has been factored out.
The argument of the exponential in the integrand is
\begin{equation}
\Phi (\tau)=(\cosh \tau-1)\cos\sigma+2\sin\left(\Frac{\theta-\sigma}
{2}\right)\sin\left(\Frac{\theta +\sigma}{2}\right)+ (\sigma
-\theta)\sin\theta \,.
\end{equation}
\noindent This formula corrects a misprint in \cite{Gil1} (Eq. 33).
 The difference $\theta-\sigma$ can
be computed in a stable way for small values of $\tau$ by using the
expression. 

\begin{equation}
\sin(\theta-\sigma)=\Frac{\sin\theta}{\cos\theta\frac{\tau}{\sinh\tau}+\cos\sigma}\,
\left[1-\frac{\tau^2}{\sinh^2\tau}\right]\,,
\end{equation}
\noindent 
together with the definition of $\sigma $ (\ref{definiciones}) 
and specific algorithms to compute
$\cosh (\tau)-1$ and $1-\sinh(\tau)^2/\tau^2$ for small $\tau$. 

The non-oscillating integral representations for 
$L_{ia}(x)$ and its derivative can be written after factoring the dominant
exponential contribution as:

\begin{equation}
L_{ia} (x)=\Frac{e^{\lambda}}{2\pi}\left[\displaystyle\int_{-\theta-\pi}^{-\theta+\pi}
e^{x\gamma (\sigma )}d\sigma-(1-e^{-2\pi a})e^{-\eta}\displaystyle\int_0^{+\infty}e^{-x\Phi (\tau )}\Frac{d\sigma}{d\tau}d\tau\right]
\label{lmon}
\end{equation}

\noindent where 
\begin{equation}
\begin{array}{l}
\gamma (\sigma)=2\sin\Frac{\theta -\sigma}{2}
\sin\Frac{\theta +\sigma}{2}+(\sigma-\theta) \sin\theta \\
\\
\eta=2x[\cos\theta + (\theta -\pi /2)\sin\theta]=2x\left(\sqrt{1-(a/x)^2}-\Frac{a}{x}\arccos(\Frac{a}{x})\right)
\end{array}
\end{equation}
 and using 
Eq. (\ref{definiciones}) 
\begin{equation}
\Frac{d\sigma}{d\tau}=\tan \sigma \left[
\Frac{1}{\tau}-\coth \tau\right] .
\end{equation}

 The first integral is dominant over
the second one for large values of the parameters and $a/x$ not too close to $a=x$. As $a\rightarrow x$ both
integrals become of the same order. 

Similarly, we have the following representation for $L'_{ia}(x)$:

\begin{equation}
L'_{ia} (x)=\Frac{e^{\lambda}}{2\pi}\left[\displaystyle\int_{-\theta-\pi}^{-\theta+\pi}
\cos\sigma e^{x\gamma (\sigma )}d\sigma+(1-e^{-2\pi a})e^{-\eta}\displaystyle\int_0^{+\infty}e^{-x\Phi (\tau )}h(\tau) d\tau\right]
\label{lpmon}
\end{equation}
\noindent

\noindent where
\begin{equation}
h(\tau)=\sin\sigma \left[\Frac{\cosh \tau}{\tau}-\Frac{1}{\sinh\tau}\right]
\end{equation}

These integral representations for $L_{ia}(x)$ and $L'_{ia}(x)$ can be used for checking the 
computation of these functions in the monotonic region. They are not used by our
algorithms \cite{Gil0} because the Airy-type asymptotic expansion (Section \ref{atuae}) and the 
expansion for large $x$ (Section \ref{asyx})
are sufficiently accurate for this functions and they are faster to compute (see \cite{Gil0}, Section 2).

\subsubsection{Oscillatory case ($x<a$)}

The non-oscillating integral representations for the oscillatory region
 are more 
difficult to evaluate numerically than those for the monotonic case. 
Indeed, as it was discussed in \cite{Gil1}, the steepest descent method
leads to three integrals, which have to be computed separately. However,
as we later discuss, for moderately large $a$ it will be enough to
compute a single integral.

In \cite{Gil1}, the following formula was obtained:

\begin{equation} 
\begin{array}{ll} K_{ia}(x)&=e^{-\pi a/2}  
\left[\displaystyle\int_{\mu}^{\infty} e^{-\Psi(\tau)} 
\left(  \cos\chi +\sin\chi\Frac{d\sigma}{d\tau}      \right)d\tau \right.\\
& \\
 & +\Frac{1}{\sinh \pi a} \displaystyle\int_{\mu-\tanh
\mu}^{\mu} \left(\cos\chi\, \sinh\rho +
\sin\chi\,\cosh\rho\,\Frac{d\sigma}{d\tau} \right)d\tau \\
\\ &\left. -\Frac{1}{\sinh \pi a} \displaystyle\int_{\pi}^{3\pi /2} \left(\cos\chi\,
\sinh\rho \,\Frac{d\tau}{d\sigma} + \sin\chi\,\cosh\rho
\right)d\sigma\right]\,,
\end{array} 
\label{oscil}
\end{equation}

where $\chi = x \sinh\,\mu-a \mu$,  
$\cosh\,\mu =\Frac{a}{x},\ \mu>0$,

\begin{equation} \Psi(\tau)=x \cosh\,\tau\,\cos\,\sigma
+a\left(\sigma-\frac12\pi\right), \rho (\tau)=-\Psi(\tau)
+a\pi\end{equation}

and

\begin{equation} \sin\,\sigma=\Frac{(\tau-\mu)\cosh\,\mu
+\sinh\,\mu}{\sinh\,\tau}\,. \label{sigma} \end{equation}

Notice that each of the three integrals  in Eq. (\ref{oscil}) can in principle be integrated with
respect to any of the two variables $\sigma $ and $\tau$, taking into account Eq. (\ref{sigma}) together
with the fact that the integration path $\tau (\sigma)$ is such that $\tau (0)=+\infty$, $\tau (\pi/2)=\mu$,
$\tau (\pi)=\mu-\tanh \mu >\tau (3\pi /2)$; however, as discussed in \cite{Gil1} there
are strong numerical reasons for the selections made. In particular, the third integral is performed with
respect to $\sigma$ (which requires numerical inversion of (\ref{sigma})) to avoid the singularity of 
$d\sigma /d\tau$ at $\tau (3\pi /2)$. As explained in \cite{Gil1} the numerical inversion of (\ref{sigma})
in the interval $\sigma\in [\pi ,3\pi /2]$ can be efficiently performed in parallel with the numerical
integration.

Similar integral representations exist for $K_{ia}' (x)$, $L_{ia}(x)$ and $L'_{ia}(x)$. We have:

\begin{equation} 
\begin{array}{ll} K'_{ia}(x)&=e^{-\pi a/2}  
\left[\displaystyle\int_{\mu}^{\infty} e^{-\Psi(\tau)} 
\left(\cos \chi\, A(\tau)+\sin\chi\, C(\tau) \right)d\tau \right.\\
& \\
 & +\Frac{1}{\sinh \pi a} \displaystyle\int_{\mu-\tanh
\mu}^{\mu} \left(
\cos\chi \cosh \rho\, A(\tau)+\sin\chi \sinh\rho\, C(\tau)  \right)d\tau \\
\\ &\left. -\Frac{1}{\sinh \pi a} \displaystyle\int_{\pi}^{3\pi /2} \left(
\cos\chi \cosh \rho\, B(\tau (\sigma))+\sin\chi \sinh\rho\, D(\tau (\sigma ))
\right)d\sigma\right]
\end{array} 
\label{oscil2}
\end{equation}

\noindent where
\begin{equation}
\begin{array}{l}
A(\tau)=-\cosh\,\tau\cos\,\sigma + \sinh
\,\tau \sin \,\sigma \Frac{d\sigma}{d\tau}\,,\,\,B(\tau)=A(\tau)
\Frac{d\tau}{d\sigma}\,,\\ \\ C(\tau)=-\sinh \,\tau \sin \,\sigma
-\cosh \,\tau \cos \,\sigma \Frac{d\sigma}{d\tau}\,,\,\,
D(\tau)=C(\tau)\Frac{d\tau}{d\sigma}\,.
\end{array}
\end{equation}

In addition, integral representations 
for $L_{ia}(x)$ and 
its derivative are:

\begin{equation} 
\begin{array}{ll} L_{ia}(x)&=\Frac{e^{\pi a/2}}{\pi}  
\left[\Frac{1-e^{-2\pi a}}{2}\displaystyle\int_{\mu}^{\infty} e^{-\Psi(\tau)} 
\left(  \sin\chi -\cos\chi\Frac{d\sigma}{d\tau}      \right)d\tau \right.\\
& \\
 & +e^{-\pi a} \displaystyle\int_{\mu-\tanh
\mu}^{\mu} \left(\sin\chi\, \sinh\rho -
\cos\chi\,\cosh\rho\,\Frac{d\sigma}{d\tau} \right)d\tau \\
\\ &\left. -e^{-\pi a} \displaystyle\int_{\pi}^{3\pi /2} \left(\sin\chi\,
\sinh\rho \,\Frac{d\tau}{d\sigma} - \cos\chi\,\cosh\rho
\right)d\sigma\right]
\label{oscil3}
\end{array}
\end{equation}

\noindent and

\begin{equation} 
\begin{array}{ll} L^{\prime}_{ia}(x)&=\Frac{e^{\pi a/2}}{\pi}  
\left[\Frac{1-e^{-2\pi a}}{2}\displaystyle\int_{\mu}^{\infty} e^{-\Psi(\tau)} 
\left( \sin\chi\, A(\tau) - \cos\chi\, C(\tau)\right)d\tau \right.\\
& \\
 & +e^{-\pi a} \displaystyle\int_{\mu-\tanh
\mu}^{\mu} \left( \sin\chi\cosh\rho\, A(\tau) - \cos\chi\sinh\rho\, C(\tau) \right)d\tau \\
\\ &\left. -e^{-\pi a} \displaystyle\int_{\pi}^{3\pi /2} \left(
\sin\chi \cosh\rho\, B(\tau (\sigma))-\cos\chi\sinh\rho\, D(\tau (\sigma)) 
\right)d\sigma\right]
\label{oscil4}
\end{array}
\end{equation}

Notice that the dominant exponential behaviour has been factored for both the
functions $K_{ia}(x)$ and $L_{ia}(x)$ and their derivatives, which
coincides with the exponential behaviour of the uniform asymptotic expansion.
This is an interesting feature when computing scaled functions in order
to avoid overflows and/or underflows in the computation. After factoring the
dominant exponential terms ($e^{\pm \pi a /2}$), the overflow and/or underflow
problems are eliminated; notice, however,
 that when computing the integrals over
finite intervals we should evaluate $\sinh(\rho)/e^{a\pi}$, 
$\cosh(\rho)/e^{a\pi}$ for $L_{ia}(x)$ and its derivative and 
$\sinh(\rho)/\sinh{a\pi}$, 
$\cosh(\rho)/\sinh{a\pi}$ for $K_{ia}(x)$ (and $K_{ia}' (x)$)
 instead of computing the hyperbolic and
the exponential separately (otherwise, overflows will take place for
moderately large $a$). For this reason it is convenient to use the expressions
\begin{equation}
\begin{array}{ll}
\Frac{\cosh\rho}{\sinh \pi a}=e^{-\psi}\Frac{1+e^{-2\rho}}
{1-e^{-2\pi a}}, &
\Frac{\sinh\rho}{\sinh \pi a}=e^{-\psi}\Frac{1-e^{-2\rho}}
{1-e^{-2\pi a}}
\end{array}
\end{equation}
\noindent in the evaluation of Eqs. (\ref{oscil}) and (\ref{oscil2})
and to proceed in the same way for 
$e^{-\pi a} \cosh\rho$ and 
$e^{-\pi a} \sinh\rho$ in Eqs. (\ref{oscil3}) and (\ref{oscil4}). 
Notice that in the oscillatory region
$\rho>0$ and that for large $a$ and $x$ both $e^{-2\rho}$ and 
$e^{-2\pi a}$ will underflow. These underflow problems 
can be easily avoided by neglecting these exponential terms for large 
parameters. 

In addition, when both exponentials become negligible, the 
integral over sigma becomes negligible and the remaining two integrals
can be approximated by only one integral. We can write

\begin{equation}
\begin{array}{l}
 K_{ia}(x)\approx e^{-\pi a/2}\left[
\displaystyle\int_{\tau_0}^{\infty} e^{-\Psi(\tau)} 
\left(  \cos\chi +\sin\chi\Frac{d\sigma}{d\tau}      \right)d\tau
+{\cal O}(e^{-\pi a/2})\right] 
 \\
\\
K'_{ia}(x)\approx e^{-\pi a/2}\left[  
\displaystyle\int_{\tau_0}^{\infty} e^{-\Psi(\tau)} 
\left(\cos \chi\, A(\tau)+\sin\chi\, C(\tau) \right)d\tau
+{\cal O}(e^{-\pi a/2}) \right]
\\
\\
 L_{ia}(x)\approx\Frac{e^{\pi a/2}}{2\pi}\left[  
\displaystyle\int_{\tau_0}^{\infty} e^{-\Psi(\tau)} 
\left(  \sin\chi -\cos\chi\Frac{d\sigma}{d\tau}      \right)d\tau 
+{\cal O}(e^{-\pi a/2} )
\right]
\\
\\
L^{\prime}_{ia}(x)\approx \Frac{e^{\pi a/2}}{2\pi} \left[ 
\displaystyle\int_{\tau_0}^{\infty} e^{-\Psi(\tau)} 
\left( \sin\chi\, A(\tau) - \cos\chi\, C(\tau)\right)d\tau
+{\cal O}(e^{-\pi a/2} )\right]
\end{array}
\label{simplificado}
\end{equation}
\noindent
where $$\tau_0=\mu-\tanh \mu .$$ 

These approximations can be used for moderately large $a$, which is the
region where integrals for the oscillatory case are employed in the code \cite{Gil0}. 

It is however useful to have the complete expressions for testing
the rest of the methods. The computation through quadrature using Eqs.
(\ref{kmon}), (\ref{lmon}), (\ref{lpmon}), (\ref{oscil}), (\ref{oscil2}), (\ref{oscil3}) and
(\ref{oscil4}),
provides a way for computing the functions in the whole $(x,a)$ plane, 
except close to $a=x$, where the integrands become non-smooth. For this reason,
they have been used for checking the algorithm, although in the
oscillatory region only Eq. 
(\ref{simplificado}) is necessary when building the numerical algorithm \cite{Gil0}.  

\subsubsection{Quadrature rule}

As reported in Goodwin \cite{Goo}, the trapezoidal rule
is a very efficient method of computation of integrals 
$\int_{-\infty}^{+\infty}f(x) dx$ for rapidly decaying integrands  $f(x)$; in
particular, it is know that the error decays as $\exp(-(\pi /h)^2)$ for integrals
of the type $\int_{-\infty}^{+\infty}e^{-x^2}g(x) dx$ with $g$ analytic in
$\{z\in\mathbb{C} : |\Im z| <\pi /h\}$. 
After appropriate changes of variable, similar arguments follow for 
integrals over finite intervals with a smooth integrand 
\cite{Sch,Tak1,Tak2}.

The semi-infinite integrals in this Section are  appropriate 
for their computation by using the trapezoidal rule, because they
 decay as a double exponential as $\tau\rightarrow
+\infty$. On the other hand, the integrals over finite intervals 
show abrupt variations as $a\rightarrow x$, 
particularly in the oscillatory case, but under an adequate change
of variables they can be also computed efficiently 
by means of the trapezoidal rule. For finite integrals, we 
consider a change of variable in order to map the finite interval
$[a,b]$ into $(-\infty ,+\infty)$
 and a successive change to improve
the convergence of the trapezoidal rule \cite{Tak1,Tak2}; namely, we
consider the following transformation:

\begin{equation}
\begin{array}{l}
I=\displaystyle\int_a^b f(x)dx=\displaystyle\int_{-\infty}^{+\infty}g(t)dt
\,,\,\,g(t)=f(x(t))
\Frac{(b-a)\cosh t}{2 \cosh^2 (\sinh t)} ,\\
\\
x(t)=\Frac{b+a}{2}+\Frac{b-a}{2}\tanh (\sinh t)\,.
\end{array}
\end{equation}

And the integral is discretized by means of the trapezoidal rule with equal
mesh size:
\begin{equation}
\displaystyle\int_{-\infty}^{+\infty}g(t)dt =h\sum_{n=-\infty}^{+\infty}
g(nh)+\epsilon \,,
\label{dis}
\end{equation}

\noindent where the error $\epsilon$ is expected to decay very fast as
the mesh size is decreased because the integrand is analytic and
its decay is doubly exponential.
 We use a trapezoidal rule which halves
the mesh size until the prescribed precision is reached; the
same rule controls that the truncation of the infinite series (\ref{dis}) gives
an error well below the accuracy claim.

Regarding the semi-infinite integrals, we use a change of variable to 
transform the integration interval $[a,+\infty)$ to $(-\infty, \infty)$. 
We consider the following change of variables to perform this map.
\begin{equation}
\tau (x)=a+\sinh^{-1}(e^{x})\,.
\end{equation}
The additional change $x=\sinh t$ improves the convergence of the
trapezoidal rule.

It is observed that, typically, no more than $8$ iterations
of the trapezoidal rule are needed, which means that the integrands are 
evaluated at $2^8 +1=257$ points in the worst cases. This is the typical 
number of
iterations for the evaluation of $\int_{-\infty}^{+\infty} e^{-x^2}dx$ 
by means of a recursive trapezoidal rule when
double precision accuracy is demanded. This fact confirms
 that the above mentioned
changes of variable are adequate for the computation of the integrals
for the modified Bessel functions.

\subsection{Continued fraction method}

 As discussed in \cite{Gil1} both $K_{ia}(x)$ and $K^{\prime}_{ia}$ can
 be computed for moderate $a$ by means of a continued fraction method,
 similar to the corresponding method for Bessel functions of real orders
 (see \cite{Tem1}  and \cite{NR}, pp. 239-240). 
We refer to \cite{Gil1} for a full description of this scheme.
 
 As numerical experiments show, this method is competitive in speed
 with asymptotic expansions  for large $x$ (Section \ref{asyx})
and the range of
application is larger. Therefore, the continued fraction method substitutes
the use of asymptotic expansions for large $x$.

 \section{Range of computation and scaled functions}
\label{range}

 As described in previous sections, the integral representations which were
developed indicate that the dominant behaviour for the functions when the 
parameters are large is of exponential type. 
This means that the computations
can only be carried for not too large values of $a$ and $x$ in order to avoid
overflows/underflows in the computation. For instance, from Eqs. (\ref{oscil})
and (\ref{oscil2}) it is seen that for large $a$ ($a>x$), we have
$K_{ia}(x)\sim e^{-a\pi /2}$ and similarly for the derivative, while for $L_{ia}(x)$
and its derivative (Eqs. (\ref{oscil3}) and (\ref{oscil4})) the asymptotic behaviour is $\sim e^{a\pi /2}$. 
This means that to avoid overflow/underflows in the computation, we must restrict the range of $a$ in the oscillatory
region to
\begin{equation}
a<2\ln (10^n N)/\pi
\end{equation}
\noindent where $N$ is either the inverse of the underflow number (when computing $K_{ia}(x)$ or its derivative) or
the overflow number (for $L_{ia}(x)$ and its derivative); 
$10^n$ is a safety factor (in the program, we take $n=8$). For processors using the IEEE standard in double precision
this will approximately limit $a$ to $a<440$. On the other hand, for the monotonic region ($x>a$) the integral
 representations
show that the dominant exponential behaviour is $K_{ia}(x)\sim e^{-\lambda}$, $L_{ia}(x)\sim e^{+\lambda}$ where 
$\lambda (x,a)=x(\cos\theta+\theta \sin\theta)$, $\sin\theta =a/x$
($\theta\in [0,\pi /2]$),
and similarly for the
derivatives. This means that, in order to avoid overflows/underflows, the range of computation must be restricted to:
\begin{equation}
\lambda (a,x)=\sqrt{x^2-a^2}+a \,\arcsin(a/x)<\ln (10^n N)\,.
\end{equation}

Figure 1 shows the computable range for $10^n N=10^{300}$ (typical value
for IEEE standard double precision)

Given that all our expressions
have the dominant exponential contributions factored out, exponentially scaled 
functions can be defined which are computable in wider ranges. Namely, we
define:
\begin{equation}
\label{escaladas}
\begin{array}{ll}
\widetilde{K_{ia}}(x)=\left\{
\begin{array}{ll}
e^{\lambda (x,a)} K_{ia}(x) & x\ge a\\
\\
e^{a\pi /2}  K_{ia}(x) & x<a 
\end{array}
\right.
&
\widetilde{K^{\prime}_{ia}}(x)=\left\{
\begin{array}{ll}
e^{\lambda (x,a)} K^{\prime}_{ia}(x) & x\ge a\\
\\
e^{a\pi /2}  K^{\prime}_{ia}(x) & x<a 
\end{array}
\right.
\end{array}
\end{equation}
\noindent
 and
\begin{equation}
\label{escaladas2}
\begin{array}{ll}
\widetilde{L_{ia}}(x)=\left\{
\begin{array}{ll}
e^{-\lambda (x,a)} L_{ia}(x) & x\ge a\\
\\
e^{-a\pi /2}  L_{ia}(x) & x<a 
\end{array}
\right.
&
\widetilde{L^{\prime}_{ia}}(x)=\left\{
\begin{array}{ll}
e^{-\lambda (x,a)} L^{\prime}_{ia}(x) & x\ge a\\
\\
e^{-a\pi /2}  L^{\prime}_{ia}(x) & x<a 
\end{array}
\right.
\end{array}
\end{equation}

\vspace*{0.3cm}
 \centerline{\protect\hbox{\psfig{file=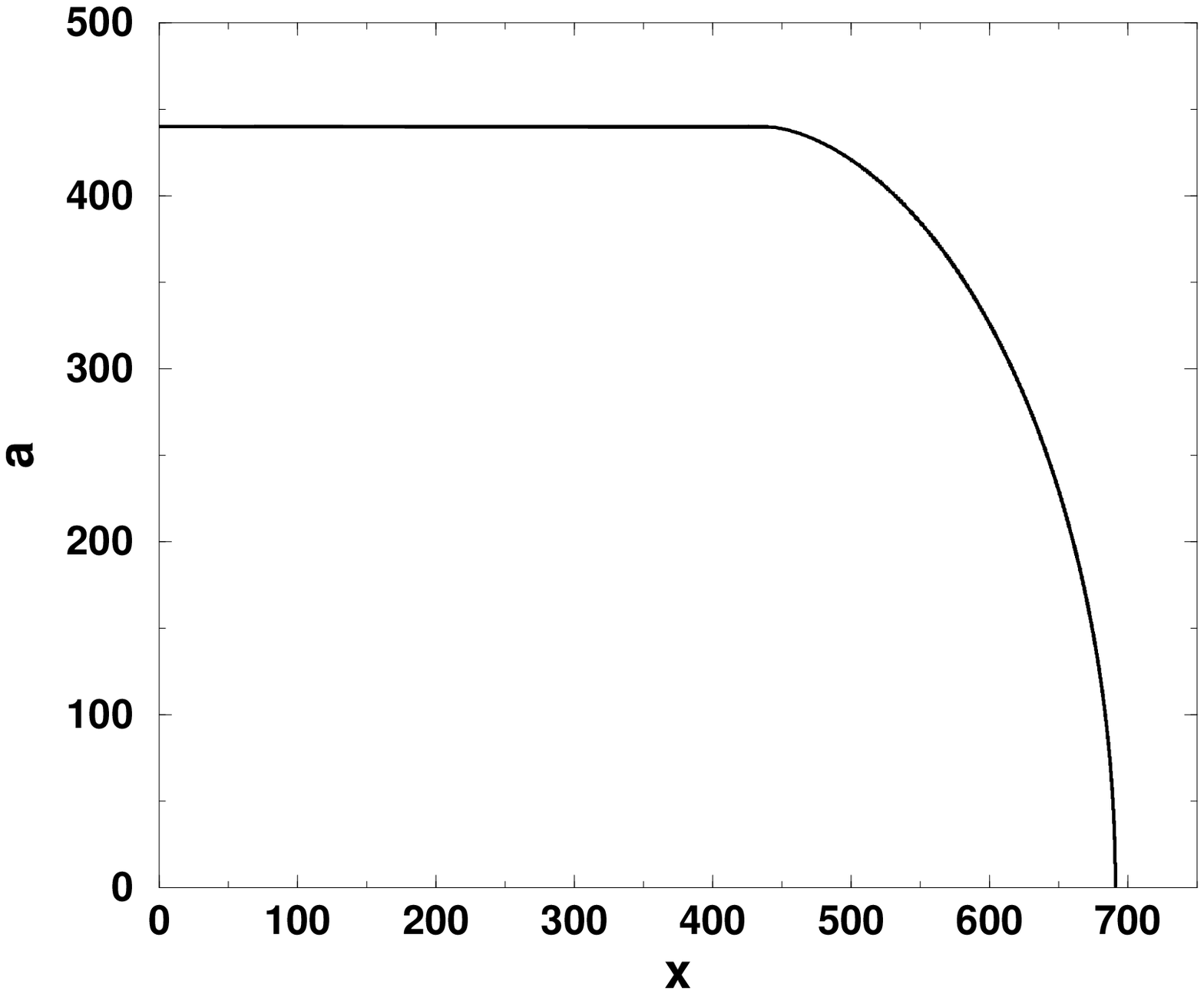,width=7cm}}}
 \vspace*{0.2cm}

 {\bf Figure 1.} {\footnotesize Computable range for the evaluation of
 $K_{ia}(x)$, $L_{ia}(x)$ and their derivatives.}

\vspace*{0.2cm}

 Note that, as in the rest of the article, we
 are considering  positive $a$ because $K_{ia}(x)$ and $L_{ia}(x)$ are even functions of $a$. 
Of course, when applying the scaling factors for negative $a$, we should replace $a$ 
by $|a|$ in the exponential scaling factors of Eqs. (\ref{escaladas}) and (\ref{escaladas2}). 
In this way, the scaled functions are also even functions of $a$. 

The definitions in (\ref{escaladas}) and (\ref{escaladas2}) eliminate exactly the front 
exponential factor in the
oscillatory region $(a>x$) from 
the series and the Airy type asymptotic expansion and in all the 
$(x,a)$ plane for the the integral 
representations. In other cases, there remains an exponential factor
with soft variation. For example, when using Airy-type
expansions in the monotonic region (neglecting non-exponential factors),
we have
\begin{equation}
\widetilde{K_{ia}}(x)\sim e^{\lambda -a\pi /2}=e^{\widetilde{\lambda}}
\end{equation}   
\noindent
where
$$\widetilde{\lambda}=x(\cos\theta +(\theta -\pi/2) \sin\theta)=\sqrt{x^2-a^2}+a (\arcsin(a/x)-\pi /2 ),$$ 
\noindent 
which is small
for $x\simeq a$ ($\theta\simeq \pi/2$); loss of accuracy in the computation of 
$\widetilde{\lambda}$ for $x\simeq a$ 
can be reduced by expanding $\widetilde{\lambda}$ in powers of $\theta-\pi /2$.

Similarly, an exponential factor remains when rescaling the 
asymptotic expansions
and the same happens when applying the continued fraction method. In this case, we 
have for $x>a$;

\begin{equation}
\widetilde{K_{ia}}(x)\sim e^{x-\lambda}=e^{-\bar{\lambda}}
\end{equation} 
\noindent
where 
$$\bar{\lambda}=x((\cos\theta -1)+\theta\sin\theta)=x-\sqrt{x^2-a^2}-a \,\arcsin(a/x),$$ 

\noindent which goes to zero as $a/x\rightarrow 0$ ($\theta\rightarrow
0$). Loss of accuracy in the computation of $\bar{\lambda }$ for small $\theta$ ($a/x$ small)
can be avoided by expanding $\bar{\lambda}$ in powers of $\theta$.

 \begin{acks}
 A. Gil acknowledges financial support from   Ministerio de 
 Ciencia y Tecnolog\'{\i}a                        
 (Programa Ram\'on y Cajal). 
 A. Gil and J. Segura acknowledge CWI 
 Amsterdam for the hospitality and financial support.  
 \end{acks}

 
 %
 %
 \end{document}